# Long-Term Voltage Stability Assessment of an Integrated Transmission Distribution System

Ankit Singhal, Venkataramana Ajjarapu
Department of Electrical and Computer Engineering
Iowa State University, Ames, IA, USA
ankit@iastate.edu, vajjap@iastate.edu

*Abstract*— **Long-term voltage stability assessment (VSA) of independent transmission (T-VSA) and distribution (D-VSA) systems have been studied since long to estimate load margins. However, their impacts on each other have been neglected due to simplified assumptions i.e. in transmission systems, loads are assumed to be aggregated, and in distribution systems, substation bus voltage is assumed to be constant. This work investigates the VSA of integrated transmission-distribution (TD-VSA) using PV curve superimposition approach and reveals the possibility that the overall system loadability may be limited by the distribution system rather than the transmission system. In this paper, we analyze why T-VSA and D-VSA are not reliable enough to estimate true load margin of the overall system. The analysis has been verified on an integrated test system in different scenarios with and without DER penetration. Overall, the paper points out the need for integrated analysis and further builds a case why it is essential to develop a realistic co-simulation framework for a reliable long-term VSA of large-scale coupled T-D systems.**

## I. INTRODUCTION

Voltage collapse phenomenon of transmission system has always been a topic of interest for power system operators and researchers. Continuously increasing load demand is forcing grid to operate at closer to the loadability limit than ever. Several power grid blackouts in past have motivated various transmission voltage stability assessment (T-VSA) studies to accurately detect voltage collapse point and estimate load margin [1], [2]. Continuation power flow (CPF) is a widely accepted method to identify precise load margin by tracing accurate PV curves [3]. Several voltage collapse indices have also been proposed to enable precise online monitoring and prediction of collapse points [4]. However, in T-VSA studies, loads have always been modeled as an aggregated load. Whereas in real-life, the loads are spread throughout the distribution systems which are connected to load buses of transmission systems. Due to aggregation of loads rather than considering full distribution systems, the impacts of electrical distances of loads have been ignored in T-VSA studies which is a potential source of error in load margin estimation.

On the other hand, voltage collapse in distribution feeders also has been identified as a critical issue since long [5]. A major blackout (June 1997) in the S/SE Brazilian system is attributed to voltage instability problem in one of the distribution networks which readily spread to the transmission grid [6]. Moreover, distribution voltage stability assessment (D-VSA) has gained significant attention recently with the arrival of distributed energy resources (DERs) on the feeders e.g. a continuation distribution power flow tool has been developed with DER integration [7]. Several attempts have been made to assess how DER penetration improves load margin of distribution systems [8]–[12]. However, in D-VSA, the substation bus has always been considered as a slack bus with an assumed constant voltage. Whereas in real-life, this assumption doesn't hold true, and the substation voltage keeps changing based on the transmission power flow. This assumption might lead to significant errors in D-VSA results.

Thus, the T-VSA and D-VSA have been studied separately as independent systems neglecting their impact on each other so far. In reality, both transmission and distribution systems are coupled physically and affect each other. Especially, at high load operating condition near the point of collapse (PoC), the substation voltage is significantly lower than the assumed constant value, thus leading to considerable exaggeration in load margin estimation of distribution feeders. Similarly, it is hard to represent distribution network losses, its power transfer limit and DER's voltage supporting capability in T-VSA. Therefore, it is essential to consider an integrated transmission-distribution (TD) system and conduct TD-VSA studies for accurate estimation of load margin of the overall system. A recent study on integrated TD system indicates that the T-VSA and D-VSA may significantly overestimate or underestimate the actual load margin with DERs [13]. However, the study does not discuss the impact of loadability limit of the distribution network in T-VSA which can not be captured by load aggregation.

Extending [13], a systematic analysis of the integrated system is presented in this paper. In this work, we propose PV curve and hypersurface superimposition approach to explain integrated TD-VSA using independent system analysis i.e. T-VSA and D-VSA. In this approach, we also discuss how the impact of substation bus voltage can be considered by taking it as another parameter along with the load increase direction and developing a two-parameter surface (or hypersurface) rather

This work was supported by the Department of Energy's Sunshot initiative program DE-0006341. The authors gratefully acknowledge the supporters' contribution.

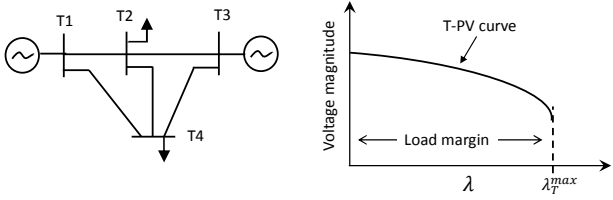

Fig. 1. A typical transmission system and PV curve

than just PV curve while evaluating D-VSA. The main contribution of the work lies in the superimposition analysis which reveals a crucial information i.e. which of the transmission and distribution systems is reaching its power transfer limit first and causing instability. The objectives of the present work are two folds: 1.) To point out the possibility of system voltage collapse being caused by the distribution system rather than the transmission system. This information can further be interpreted and used for various applications such as to take corrective actions in stressed operating conditions or for planning purpose, and 2.) To emphasize the need of integrated TD VSA to estimate true load margin of the system, and to motivate the further development of a co-simulation framework for a reliable voltage stability assessment of realistic large-scale coupled T-D systems.

The remainder of the paper is organized as follows. Section II analyzes VSA of independent transmission and distribution systems. Building on that, the TD-VSA of an integrated system and PV curve superimposition analysis are presented in Section III. Then, the analysis has been verified via numerical tests in Section IV. At the end, discussion and the concluding remarks are provided in the Section V and Section VI respectively.

## II. INDEPENDENT SYSTEMS ANALYSIS

In this section, we will analyze VSA of independent transmission (T-VSA) and independent distribution systems (D-VSA) to further develop an integrated T-D system analysis (TD-VSA) in the next section.

### A. Transmission System VSA (T-VSA)

To analyze VSA of a transmission system, parameterized transmission power flow (TPF) equations can be written by introducing the loading parameter $\lambda$ as follows:

$$G_T(x_T, \lambda) = S_T(\lambda) - S_{TT}(x_T) = 0 \qquad (1)$$

Where, $S_T$ is a vector of net complex power injections (generation-load) at all transmission system buses; $S_{TT}$ is a vector of total complex power flowing out of each bus. Loading parameter is denoted by $\lambda$ and voltage magnitude and angle variables are denoted by $x_T = [v_T, \theta_T]$. Using CPF [3], PV curve for a transmission system can be traced to assess loadability limit at $\lambda = \lambda_T^{max}$ and consequently load margin from the operating point. A typical transmission system and PV curve are shown in Fig. 1.

### B. Distribution System VSA (D-VSA)

A typical distribution system is shown in Fig. 2. Similar to (1), parameterized distribution power flow (DPF) equations can be written as follows:

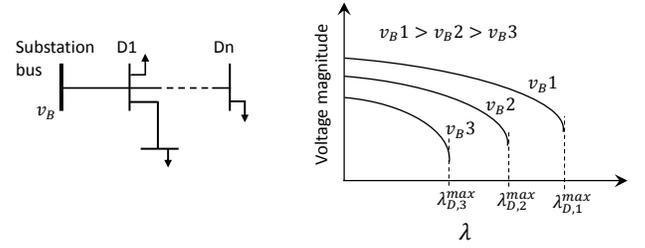

Fig. 2. A typical distribution system and PV curve on various substation voltage values

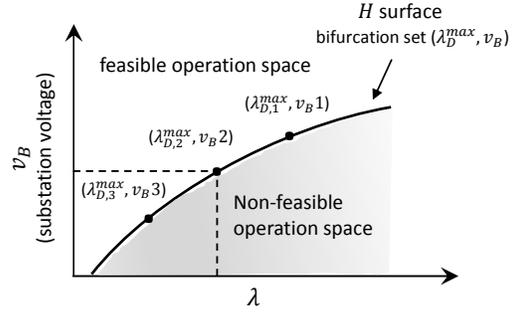

Fig. 3. An example of a hypersurface $H$ of a distribution system

$$G_D(x_D, v_B, \lambda) = S_D(\lambda) - S_{DD}(x_D, v_B) = 0 \qquad (2)$$

Notations are similar to TPF in (1) where subscript $D$ denotes the distribution system buses. It is crucial to note there are two independent parameters in DPF equation (2) i.e. loading parameter ($\lambda$) and substation bus voltage ($v_B$). For distribution system analysis, substation bus is considered a slack bus and $v_B$ is treated as a fixed parameter instead of a variable. But in real-world integrated T-D system, $v_B$ is not fixed and decided by the transmission power flow solution. Therefore, to observe the impact of $v_B$ on D-VSA, PV curve for a distribution system can be traced for different $v_B$ values as shown in Fig. 2. Higher value of substation voltage leads to higher loadability limit ($\lambda_D^{max}$) and consequently higher load margin. In other words, the distribution system has two parameters where critical value of loading parameter ($\lambda_D^{max}$) is a function of another parameter $v_B$. In a multi-parameter space, we can construct a hypersurface $H$ which is a boundary of the feasible region of operation at stable equilibrium [14]. $H$ surface is a set of substation voltage magnitude and corresponding critical loading.

A typical $H$ surface is shown in Fig. 3. Horizontal axis is loading parameter and the vertical axis is substation voltage magnitude. A point P2 on the surface indicates that if substation voltage is maintained at $v_B 2$, the load on distribution system can only be increased till $\lambda = \lambda_{D,2}^{max}$. To increase the loadability limit to $\lambda_{D,1}^{max}$, substation voltage must be increased to $v_B = v_B 1$. The space above the $H$ surface is a feasible solution region. As we approach towards the surface, we move closer to the collapse point. On the surface, we reach to nose points where a slight increase in load or decrease in $v_B$ may result in voltage collapse. The impact of decreasing substation voltage can easily be understood by plotting $H$ surface. We will analyze how $H$ surface affects integrated T-D VSA in the next section.

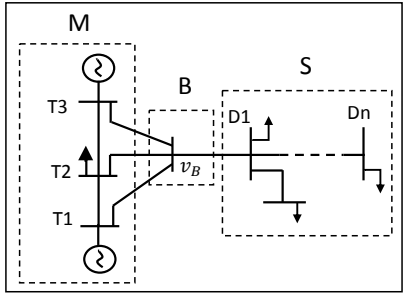

Fig. 4. A representative example of an integrated T-D system with master (M), boundary (B) and slave (S) subsystems

## III. INTEGRATED T-D SYSTEM ANALYSIS

### A. Integrated TD-VSA Formulation

To understand the role of both the transmission and distribution systems, a parameterized integrated T-D system can be modeled through master-slave splitting method [15]. In this method, the integrated T-D system is split into three subsystems as shown in a representative example in Fig. 4 i.e. master (M), boundary (B) and slave (S). M and S subsystems include all transmission and distribution systems buses, respectively excluding T-D boundary buses which are included in subsystem B. A set of parameterized power flow equations for such system can be written as follows:

$$G_{TD}(x_M, x_B, x_S, \lambda) = \begin{cases} G_M(x_M, x_B, \lambda) = 0 \\ G_B(x_M, x_B, x_S, \lambda) = 0 \\ G_S(x_B, x_S, \lambda) = 0 \end{cases} \quad (3)$$

$G_M, G_B, G_S$ are power flow equations and $x_M, x_B, x_S$ are voltage magnitude and angle variables for the master, boundary and slave subsystems respectively. To see the impact of independent T-VSA and D-VSA on integrated VSA, we can explicitly write $G_M, G_B, G_S$ as follows:

$$\begin{aligned} G_M &= S_M(\lambda) - S_{MM}(x_M) - S_{MB}(x_M, x_B) \\ G_B &= -S_{BB}(x_B) - S_{BM}(x_B, x_M) - S_{BS}(x_B, x_S) \end{aligned} \quad (4.1)$$

$$G_S = S_S(\lambda) - S_{SS}(x_S) - S_{SB}(x_S, x_B) \quad (4.2)$$

$S_X$ denotes a vector of net complex power injection of buses of system $X$; and $S_{XY}$ denotes the complex power flow from buses of system $X$ to system $Y$. In (4), we can see the components of T-VSA and D-VSA in (4.1) and (4.2) respectively. Equation (4.1) is similar to T-VSA (1) except an additional term $S_{BS}(x_B, x_S)$ which denotes the net power transfer from substation (boundary) to distribution system (slave). In T-VSA (1), this power transfer was part of $S_T(\lambda)$ in form of fixed aggregated loads. In (4.1), however, $S_{BS}$ is not an aggregated load but a separate variable being decided by (4.2). $S_{BS}$ contains information of the physical distribution network which can not be captured accurately by an aggregated load modeling such as real and reactive power losses, power transfer limitation of distribution lines etc. This leads to error in estimation by T-VSA. Similarly, (5.2) is same as D-VSA (2) except the fact that the substation voltage $v_B$ is not an independent parameter but a variable being decided by (5.1). Change in $v_B$ significantly affects the loadability limit of distribution lines which is hard to capture with the assumption of constant $v_B$ in D-VSA. This leads to error in D-VSA estimation. Thus, in an integrated T-D system, the D-VSA and T-VSA are coupled through two variables at boundary bus i.e. voltage magnitude and net power transfer from distribution system to boundary bus. A realistic capture of these two variables leads to a true estimation of load margin in TD-VSA compare to T-VSA and D-VSA. By using CPF on (4), we can trace the PV curve for T-D system and estimate loadability limit at $\lambda = \lambda_{TD}^{max}$.

### B. PV Curve Superimposition Analysis

Since, in an integrated system, voltage at the boundary bus $v_B$ is a coupling factor and loading parameter $\lambda$ is the same for both T and D system buses, we can superimpose $H$ surface of the distribution system as shown in Fig. 3 with T-PV curve of transmission system as shown in Fig. 1. Please note that the y-axis in T-PV curve is voltage at the same boundary bus $v_B$. The superimposition of the two curves reveals some useful inferences regarding T-D VSA. Superimposition provides possibility of the two cases which are explained in detail below:

*1) Case A: Constrained by Distribution System*

This is the case where distribution $H$ surface intersects T-PV curve at $\lambda_{TD}^{cut} < \lambda_T^{max}$ on superimposition as shown in Fig. 5. Independent transmission system (T-PV curve) can go till $\lambda_T^{max}$ but since the coupled TD system can not cross the $H$ surface and enter into infeasible grey region, it has to stop before the intersection $\lambda_{TD}^{cut}$. Physically, the voltage at the boundary bus goes lower than the minimum substation voltage distribution feeder can handle at that particular load level and thus the distribution loadability limit arrives before the transmission system loadability limit. So, this is the case where overall system voltage collapse is being caused by the distribution systems. It should be noted since the $S_{BS}$ includes distribution feeder losses in TD-VSA, the TD-PV curve does not exactly follow the T-PV curve. As load increase, losses increase and TD-PV curve moves away from the T-PV curve more. The analysis indicates the actual loadability limit will be less than $\lambda_{TD}^{cut}$ i.e. $\lambda_{TD}^{max} \leq \lambda_{TD}^{cut}$. If losses in distribution system are to be considered zero or negligible, then $\lambda_{TD}^{max} \approx \lambda_{TD}^{cut}$.

*2) Case B: Constrained by Transmission System*

This is the case when the distribution $H$ surface does not intersect the T-PV curve as shown in Fig. 6. Physically, even at very high loading condition, the distribution system can handle the low substation voltage in this case. Consequently, transmission system hits the loadability limit before the distribution system. It infers that the distribution feeder is not being the limiting factor in maximum loadability of the coupled system. Since, the loadability limit is caused mainly by the transmission system, $\lambda_{TD}^{max}$ is very close to $\lambda_T^{max}$ but not exactly same because of distribution losses taken into account. Generally, $\lambda_{TD}^{max} \leq \lambda_T^{max}$ but if losses can be neglected, then $\lambda_{TD}^{max} \approx \lambda_T^{max}$.

Though the exact $\lambda_{TD}^{max}$ can be estimated by solving (5), this analysis provides additional information that which system is constraining the overall coupled system so that the appropriate actions can be taken.

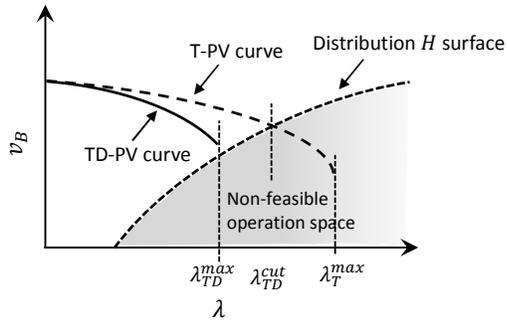

Fig. 5. Case A: distribution hypersurface intersects T-PV curve i.e. distribution loadability is limiting factor

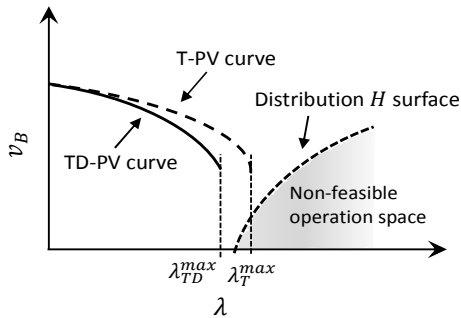

Fig. 6. Case B: distribution hypersurface doesn't intersect T-PV curve i.e. transmission loadability is limiting factor

### C. Impact of DERs on the Integrated System VSA

DERs affect the net load, thus directly impact the net power transfer from distribution feeder to the transmission system. Though, some of these impacts can be modeled in T-VSA such as net load reduction but many other impacts are very difficult to consider such as the significant impact on losses in the feeder, improvement in distribution system loadability etc. It can exaggerate or underestimate the load margins. Therefore, integrated TD-VSA becomes essential and provides a more accurate assessment. Most importantly in cases where loadability is limited by distribution systems, integrated analysis can evaluate how the DERs affect the $H$ surface of feeders and consequently affect overall system loadability. DERs have two major effect on the load margin. First, they decrease the net real power demand of the system which shifts the operating point backwards in terms of loading. This should increase the load margin by the same amount. Second, they also reduce the net PF of the substation as they only decrease real power demand while not altering the var demand when operating in unity PF mode which is a widely accepted practice. And as we know, as net substation PF decreases, PV curves shifts downwards and loadability limit decreases. This leads to decrease in load margin. Thus, the net effect of DERs is a combination of the two opposite factors. Generally, the load margin may increase but because of reduction in net substation PF, the increase would be less than the expected amount i.e. decreases in the net load.

## IV. NUMERICAL TESTS

An integrated T-D system has been constructed to verify the analysis. Load bus 5 (90 MW) of IEEE 9 bus transmission system is expanded by attaching several distribution feeders to it. 2 type of distribution systems D1 and D2 are used which are modified balanced versions of IEEE 4 bus test system [16]. To match the load demand of transmission load bus, distribution feeders are duplicated multiple times. In the results, load margins of the integrated system by TD-VSA are compared with the same obtained from T-VSA to evaluate the need for integrated assessment. CPF is used to evaluate VSA and trace PV curves by increasing loads at all the nodes with the same $\lambda$ such that the base operating case is at $\lambda = 0$.

### A. Comparison with T-VSA

Two cases of integrated systems are created i.e. case A and case B by connecting D1 and D2 respectively. Fig. 7 shows the PV curve analysis of case A. It can be seen in the Fig. 7 that the loadability limit estimated by T-VSA is $\lambda_T^{max} = 0.99$. However, in this case, the $H$ surface of the distribution system is cutting the T-PV curve at $\lambda_{TD}^{cut} = 0.59$. This represents the case where distribution system limits the overall loadability limit to $\lambda_{TD}^{max} = 0.35$. In other words, load margin of bus 5 estimated from non-integrated system is 133 MW whereas the actual load margin of the integrated system is 46 MW. In this case, T-VSA has significantly exaggerated (almost 3 times) the load margin estimation. Similar PV curve analysis of an integrated system in case B is compared with T-VSA in Fig. 8. In this case, the distribution $H$ surface does not cut the T-PV curve which infers that the loadability of the integrated system is not constrained by the distribution feeders rather by the transmission system itself. TD-VSA estimates $\lambda_{TD}^{max} = 0.88$ which means the actual load margin of 119 MW. Though, in this case also T-VSA overestimate the actual load margin but the accuracy is much better than the case A. The results are quantified in Table I. These results reveal that the actual load margin can not be evaluated correctly without considering integrated T-D impact in both the cases. These errors are mainly because of not considering distribution feeder loadability limit in case A and distribution system losses in case B. The error in case B can still be tolerated or corrected by considering a loss factor, but the error in case A presents a worst-case scenario and brings out the real necessity of the integrated analysis.

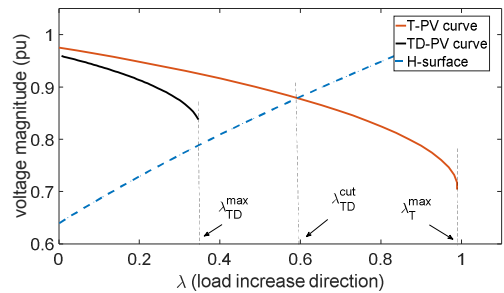

Fig. 7. PV curve assessment of an integrated system in case A where the loadability limit is constrained by the distribution system.

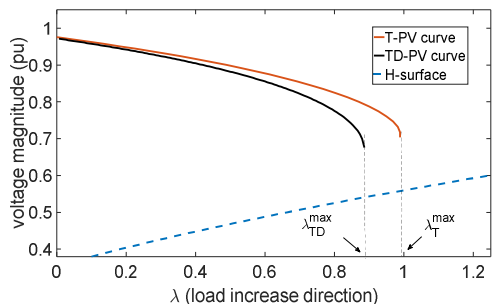

Fig. 8. PV curve assessment of an integrated system in case B where the loadability limit is constrained by the transmission system.

TABLE I
TRANSMISSION SYSTEM LOAD MARGIN (MW) COMPARISON OF BUS 5 FOR CASE A AND CASE B

| Load margin (MW) from | Case A | Case B |
|---|---|---|
| T-VSA | 133 | 133 |
| TD-VSA | 47 | 120 |

### B. Comparison with D-VSA

Integrated TD-VSA also provides a better estimation of distribution system loadability compare to D-VSA. The results are compared in Table II. In D-VSA, a constant $v_B = 1$ pu is assumed. In both the case A and case B, D-VSA exaggerate actual load margin of distribution feeder by 3 and 5 times respectively. The main reason of this overestimation in case A is assumption of constant $v_B$ in D-VSA which actually decreases in integrated TD-VSA. However, in case B, along with decreasing $v_B$ there is another reason for this error in estimation i.e. impact of transmission lines loadability. In case B, D-VSA gives a false impression of high load margin of the feeder (31.5 MW). Whereas, in reality (TD-VSA), transmission system reaches it limit much before distribution feeder can hit its limit, thus reducing the practical load margin to 6 MW. Therefore, even if the impact of decreasing $v_B$ is taken into consideration through tracing hypersurface, impact of transmission system reaching its limit first can only be captured in integrated TD-VSA.

TABLE II
DISTRIBUTION SYSTEM LOAD MARGIN (MW) COMPARISON OF D1 AND D2 FOR CASE A AND CASE B

| Load Margin (MW) from | D1 feeder (Case A) | D2 feeder (Case B) |
|---|---|---|
| D-VSA | 13 | 31.5 |
| TD-VSA | 4.8 | 6 |

### C. Impact of Load Shedding

By identifying whether the distribution or the transmission system is limiting the maximum loadability, an appropriate corrective action can be taken to improve load margin such as load shedding. Case A where distribution system is reaching its limit first should get maximum benefit from load shedding. To demonstrate that, 15 MW of load is shed at distribution system D1 and D2 in case A and case B respectively. Impact on load margin of the integrated system is compared in Table III. As expected, load margin is observed to be improved in case A because of relieving distribution feeder. In Case B, the load margin decreased rather than improving because the overall system was not limited by the distributions system on which load was shed. Thus, integrated TD-VSA provides the useful information to take corrective actions which can not be obtained by T-VSA.

TABLE III
IMPACT OF LOAD SHEDDING AT DISTRIBUTION FEEDER ON BUS 5 LOAD MARGIN (MW) IN CASE A AND CASE B

| Load margin (MW) from TD-VSA | No load shedding | Load shedding |
|---|---|---|
| Case A | 47 | 58 |
| Case B | 120 | 111 |

### D. Impact of DER

Five different scenarios with 10% - 50% DER penetration are created with both case A and case B. To understand the DER impact on VSA, TD-PV curve for case A with 50% DER penetration is compared with no DER case in Fig. 9. These are the PV curves drawn from TD-VSA analysis. The black curve is no DER case while orange curve represents DER case. As expected, there are two main effects of DER penetration as seen in Fig. 9. First, the operating point moves backward due to a reduction in net load at bus 5 due to DER generation. Second, the voltage collapse point also moves backward as the net PF at bus 5 decreases. This leads to less than expected increase in load margin due to DER e.g. in Table IV, 50% (45 MW) increase of DER in case A causes only 23 MW increase in load margin. The load margins in all cases are compared in Table IV. Another important thing to notice is that case A has much higher impact of DER than the case B. This is due to load margin being constrained by the distribution systems in case A and the DER penetration helps in relieving that constraint by increasing $\lambda_{TD}^{cut}$ of $H$-surface. This phenomenon can not be captured by T-VSA which shows almost no impact of DER at load margin (1st row).

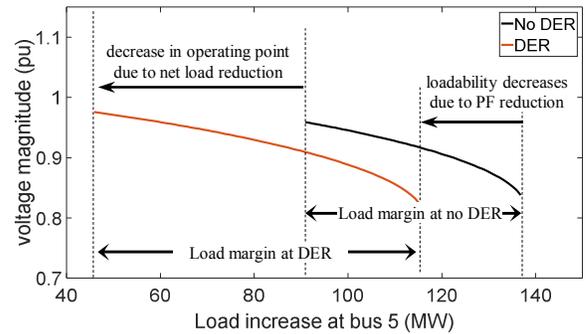

Fig. 9. Comparison of PV curve of integrated system with 50% DER penetration at both unity PF mode and VCM with no DER case

TABLE IV
BUS 5 LOAD MARGIN COMPARISON OF CASE A AND CASE B WITH DER PENETRATION AT UNITY PF

| Load margin (MW) from | DER penetration level (%) at unity PF | | | | | |
| | 0 % | 10% | 20% | 30% | 40% | 50% |
|---|---|---|---|---|---|---|
| T-VSA | 133 | 133 | 134 | 134 | 134 | 134 |
| Case A: TD-VSA | 47 | 52 | 57 | 61 | 66 | 70 |
| Case B: TD-VSA | 120 | 121 | 124 | 126 | 128 | 129 |

## V. Discussion

Results reveal that there are significant errors in loadability assessment in all scenarios if we neglect the distribution network, making the traditional T-VSA untrustworthy. Specifically, the case A comes out as the worst-case scenario (around 200% error in load margin estimation) where distribution network loadability limit constraints the overall system loadability. Moreover, the load margin improvement due to load shedding and DER penetration is significantly higher in case A, 30% and 45% respectively, compared to case B. This confirms the importance of identifying case A which can only be explained by including distribution system in the analysis.

The main intent of this work was to indicate the importance of including distribution system in VSA analysis. Taking a cue from this study, we acknowledge the need for further in-depth investigation and exploration of this phenomena on realistic large systems. Due to the large difference in impedance parameters and different nature of transmission and distribution power flow solvers, integrated power flow analysis may lead to numerical problems. Therefore, a robust co-simulation framework is more suitable rather than the integrated TD formulation to analyze a realistic system with large unbalanced 3-phase distribution systems connected to multiple load buses of transmission systems. The results from this study provoke new questions and motivate authors to further investigate the complexities of this problem via co-simulation such as unbalanced feeders, the impact of tap changers and voltage limits at distribution system, the impact of DER volt/var control etc. The study is in progress and the findings will be reported in the next paper.

## VI. Conclusion

Loadability assessment of a transmission-distribution integrated system has been analyzed in this paper and compared with the traditional T-VSA and D-VSA of independent transmission and distribution systems respectively. The superimposition analysis reveals both the possibilities of integrated TD system loadability being limited by either distribution system (case A), or the transmission system (case B). Although in both the cases, T-VSA over-estimate load margin, the error is much higher in case A because the distribution network loadability limit can not be captured in aggregated modeling in T-VSA. To confirm further, load shedding at distribution feeder and DER penetration significantly increase the TD load margin in case A because of improvement in distribution feeder loadability. Whereas it doesn't affect the load margin much in case B as the distribution loadability is not a limiting factor and improving it is not crucial for overall load margin.

Overall, the analysis and superimposition approach presented here successfully demonstrates the importance of integrated TD analysis. However, this also opens the need for a realistic co-simulation approach as discussed in Section V. New insights from this analysis can be utilized to develop online monitoring index, design corrective actions in stressed operating conditions etc. as the future scope of the work.

## Appendix

Distribution feeder impedance $(r + jx)\Omega$/mile data:

|  | Line 1-2 | Line 2-3 | Line 3-4 |
|---|---|---|---|
| Feeder D1 | 0.45+1.07j | 0.45+1.07j | 0.45+1.07j |
| Feeder D2 | 0.36+0.53j | 0.36+0.32j | 0.36+0.64j |